\documentclass{article}

\usepackage[cp850]{inputenc}
\usepackage{amssymb}
\usepackage{latexsym}

\def\Z{\mathbb{Z}}
\def\N{\mathbb{N}}
\def\R{\mathbb{R}}
\def\F{{\cal F}}
\def\P{{\cal P}}
\def\Q{{\cal Q}}
\def\S{{\cal S}}
\def\A{{\bf A}}
\def\B{{\bf B}}
\def\b{{\bf b}}

\def\t{{\bf t}}
\def\x{{\bf x}}

\begin{document}
\setlength{\parindent}{0pt}
\setlength{\parskip}{0.4cm}

\newtheorem{theorem}{Theorem}
\newtheorem{corollary}[theorem]{Corollary}
\newtheorem{definition}{Definition}

\begin{center}

\Large{\bf Multidimensional Ehrhart Reciprocity}
\footnote{Appeared in {\it J.~Combin.~Theory Ser.~A} {\bf 97}, no.~1 (2002), 187--194. \\ 
          This work is part of the author's Ph.D. thesis. \\
          {\it Mathematical Reviews Subject Numbers:} 05A15, 11D75.}

\vspace{5mm}

\normalsize

{\sc Matthias Beck}

\end{center}

\vspace{5mm}

\footnotesize
{\bf Abstract.} In \cite{beck}, the author generalized Ehrhart's idea (\cite{ehrhart}) of counting 
lattice points in dilated rational polytopes: Given a rational polytope, that is, a 
polytope with rational vertices, we use its description as the intersection of 
halfspaces, which determine the facets of the polytope. Instead of just a single 
dilation factor, we allow different dilation factors for each of these facets. We proved 
that, if our polytope is a simplex, the
lattice point counts in the interior and closure of such a {\it vector-dilated} simplex
are quasipolynomials satisfying an Ehrhart-type reciprocity law. This generalizes the classical
reciprocity law for rational polytopes (\cite{ehrhart}, \cite{macdonald}).
In the present paper we complete the picture by extending this result to general rational polytopes. 
As a corollary, we also generalize a reciprocity theorem of Stanley (\cite{stanley}). 

\normalsize


\section{Introduction}
Based on the elementary identity 
  \begin{equation}\label{gi} \left[ \frac{ t-1 }{ a }  \right] = - \left[ 
\frac{ -t }{ a }  \right] - 1 \ , \end{equation}
where $ a \in \N $, $ t \in \Z $ and $ [x] $ is the greatest integer function, we proved in \cite{beck} a generalization 
of the Ehrhart-Macdonald reciprocity law for rational polytopes. 
(A {\bf rational polytope} is a polytope whose vertices are rational.) 
More precisely, let ${\cal P}$ be an $n$-dimensional rational polytope in $ \R^{ n } $.
For a positive integer $ t $, let 
  \[ L ( {\cal P}^{\circ}  , t ) = \# \left( t {\cal P}^{\circ}  \cap \Z^{ n }  \right) \qquad \mbox{ and } \qquad L ( \overline{\cal P} , t ) = \# \left( t \overline{\cal P} \cap \Z^{ n }  \right) \] 
denote the number of integer points (''lattice points'') in the interior of the dilated polytope 
$ t {\cal P} = \{ tx : x \in {\cal P}  \} $ and its closure, respectively. 
Ehrhart, who initiated the study of the lattice point count in dilated 
polytopes (\cite{ehrhart}), proved that $ L ( {\cal P}^{\circ}  , t ) $ and 
$ L ( \overline{\cal P} , t ) $ are quasipolynomials in $t$. 
(A {\bf quasipolynomial} is an expression of the form 
  \[ c_{n}(t) \ t^{n} + \dots + c_{1}(t) \ t + c_{0}(t) \ , \] 
where $c_{0}, \dots , c_{n}$ are periodic functions in $t$.) 
He conjectured the following reciprocity law, which was first proved by Macdonald (\cite{macdonald}): 
\begin{theorem}[Ehrhart-Macdonald reciprocity law]\label{recipr} Suppose the rational po\-lytope $ \P $ is
homeomorphic to an $n$-manifold. Then 
  \[ L ( {\cal P}^{\circ}  , -t ) = (-1)^{ n } L ( \overline{\cal P} , t )  \ . \] 
\end{theorem} 
In \cite{beck}, we generalized the notion of dilated polytopes: we use the description of a convex polytope as the 
intersection of halfspaces, which determine the facets of the polytope. Instead of dilating the polytope 
by a single factor, we allow different dilation factors for each facet, such that the combinatorial type of 
the polytope does not change. Recall that two polytopes are {\bf combinatorially equivalent} if there exists 
a bijection between their faces that preserves the inclusion relation. 

It is a crucial fact that rational polytopes can be described by inequalities with integer coefficients. 
The following definition appeared in \cite{beck} only for simplices: 
\begin{definition}\label{def} Let the convex rational polytope $ \P $ be given by
   \[ \P = \left\{ \x \in \R^{ n } : \ \A \ \x \leq \b \right\} \ , \]
with $ \A \in M_{ m \times n } ( \Z ) , \b \in \Z^{ m } $. Here the inequality is
understood component\-wise. For $ \t \in \Z^{ m } $, define the {\bf vector-dilated polytope}
$ \P^{ ( \t ) }  $ as
   \[ \P^{ ( \t ) } = \left\{ \x \in \R^{ n } : \ \A \ \x \leq \t \right\} \ . \]
For those $ \t $ for which $ { P^{ ( \t ) }   } $ is combinatorially equivalent to 
$ \P = \P^{ ( \b ) } $, we define the number of lattice points in the interior and closure of
$ \P^{ ( \t ) }  $ as
   \[ i_{\P} ( \t ) = \# \left( { \P^{ ( \t ) } }^{\circ}  \cap \Z^{ n } \right) \quad \mbox{ and } \quad j_{\P} ( \t ) = \# \left( \P^{ ( \t ) } \cap \Z^{ n }  \right) \ , \]
respectively.
\end{definition} 
Geometrically, we fix for a given polytope the normal vectors to its facets and consider all
possible positions of the facets that do not change the face structure of the polytope.
Note that the dimension of $ \t $ is the number of facets of the polytope. 
The previously defined quantities $ L ( {\cal P}^{\circ}  , t ) $ and $ L ( \overline{\cal P} , t ) $
can be recovered from this new definition by choosing $ \t = t \b $. 
In \cite{beck}, we obtained a reciprocity law for vector-dilated {\it simplices}:  
\begin{theorem}\label{multithm} Let $ \S $ be an $n$-dimensional rational simplex. Then 
$ i_{\S} ( \t ) $ and $ j_{\S} ( \t ) $ are quasipolynomials in $ \t \in \Z^{ n+1 } $, satisfying 
   \[ i_{\S} ( - \t ) = (-1)^{ n } j_{\S} ( \t ) \ . \]
\hfill {} $\Box$
\end{theorem}
A {\bf quasipolynomial} in the $d$-dimensional variable $\t = ( t_{1} , \dots , t_{d} )$ is the natural 
generalization of a quasipolynomial in a 1-dimensional variable: namely, an expression of the form 
  \[ \sum_{ 0 \leq k_{1} , \dots , k_{d} \leq n } c_{ ( k_{1} , \dots , k_{d} ) } \ t_{1}^{ k_{1} } \cdots t_{d}^{ k_{d} } \ , \] 
where $ c_{ ( k_{1} , \dots , k_{d} ) } = c_{ ( k_{1} , \dots , k_{d} ) } ( t_{1} , \dots , t_{d} ) $ 
is periodic in $ t_{1} , \dots , t_{d} $. 
In \cite{beck}, we gave an actual example of such a quasipolynomial arising from a lattice point count in a polytope. 

In the present paper, we finish the picture by extending Theorem \ref{multithm} to general rational polytopes. 
We should extend Definition \ref{def} to non-convex polytopes. This can be done naturally in an additive way: 
write the polytope as the union of convex polytopes, and apply the above Definition \ref{def} to these components. 
More thoroughly, we make the following 
\begin{definition} Let $ \P $ be a rational polytope. Write $ \P = \bigcup_{ k=1 }^{ r } \P_{ k } $, where 
$ \P_{ k } $ are convex rational polytopes, say, 
  \[ \P_{ k } = \left\{ \x \in \R^{ n } : \ \A_{k} \ \x \leq \b_{k} \right\} \ , \] 
with $ \b_{ k } \in \Z^{ m_{k} } $. Given $ \t \in \Z^{ m } $, where $ m = m_{1} + \dots + m_{r} $, combine 
the first $ m_{1} $ components of $ \t $ in a vector $ \t_{1} $, the next $ m_{2} $ components in $ \t_{2} $, etc. 
Define the {\bf vector-dilated polytope} $ \P^{ ( \t ) }  $ as
   \[ \P^{ ( \t ) } = \bigcup_{ k=1 }^{ r } \P_{ k }^{ ( \t_{k} ) } \ . \] 
For those $ \t $ for which $ { P^{ ( \t ) }   } $ is combinatorially equivalent to 
$ \P $, we define as above 
   \[ i_{\P} ( \t ) = \# \left( { \P^{ ( \t ) } }^{\circ}  \cap \Z^{ n } \right) \quad \mbox{ and } \quad j_{\P} ( \t ) = \# \left( \P^{ ( \t ) } \cap \Z^{ n }  \right) \ . \] 
\end{definition} 
Finally, we derive a generalization of the following theorem of Stanley (\cite{stanley}) in terms of vector-dilated 
polytopes. 
The Ehrhart-Macdonald reciprocity law compares the lattice point count of the polytope with that of the interior, that
is, the polytope with all its facets removed. Stanley's theorem tells us what to expect if we only remove {\it some} of
the facets.
\begin{theorem}[Stanley]\label{stan} Suppose the rational polytope $ \P $ is
homeomorphic to an $n$-manifold. Denote the set of all (closed) facets of $\P$ by $F$, and let $T$ be a subset 
of $F$, such that $ \bigcup_{ \F \in T } \F $ is homeomorphic to an $(n-1)$-manifold. Let 
  \[ j_{\P, T} ( t ) = \# \left( t \left( \P - \bigcup_{ \F \in T } \F \right) \cap \Z^{ n } \right) \] 
and 
  \[ i_{\P, T} ( t ) = \# \left( t \left( \P - \bigcup_{ \F \in F-T } \F \right) \cap \Z^{ n } \right) \ . \] 
Then 
   \[ i_{\P, T} ( - t ) = (-1)^{ n } j_{\P, T} ( t ) \ . \]
\hfill {} $\Box$
\end{theorem}
Note that Theorem \ref{recipr} is the special case $ T = \emptyset $ of Theorem \ref{stan}. 
For an example that this result does not hold in general, see \cite{stanley}. 


\section{Extending Ehrhart reciprocity}
In \cite{beck}, we remarked that Theorem \ref{recipr} follows directly from Theorem \ref{multithm}. Since we will 
use Theorem \ref{recipr} to show the main result of this paper, we start by actually proving this remark. 

{\it Proof of Theorem} \ref{recipr}. We use double induction on the dimension of the polytope $\P$ and on the number of $n$-dimensional
simplices which triangulate $\P$. It is easy to see (\cite{beck}) that Theorem \ref{recipr} follows for 1-dimensional 
polytopes (that is, intervals) from (\ref{gi}). Also, Theorem \ref{recipr} holds for simplices, as a special case 
of Theorem \ref{multithm}. For a general $\P$ satisfying the hypotheses of the statement, write 
  \[ \P = \P_{1} \cup \P_{2} \ , \] 
where $ \P_{1} $ is an $n$-dimensional simplex such that $ \P_{2} := \overline{ \P - \P_{1} } $ is again a polytope homeomorphic to an 
$n$-manifold. Note that the conditions on $ \P $ imply that $ \P_{1} $ and $ \P_{2} $ share an $(n-1)$-dimensional 
polytopal boundary, which we denote by $ \P_{3} $. Hence 
  \[ L ( \overline{\P} , t ) = L ( \overline{\P_{1}} , t ) + L ( \overline{\P_{2}} , t ) - L ( \overline{\P_{3}} , t ) \] 
and 
  \[ L ( \P^{\circ} , t ) = L ( \P_{1}^{\circ} , t ) + L ( \P_{2}^{\circ} , t ) + L ( \P_{3}^{\circ} , t ) \ . \] 
By induction, we can apply Theorem \ref{recipr} to $ \P_{1} $, $ \P_{2} $, and $ \P_{3} $: 
  \begin{eqnarray*} &\mbox{}& L ( \overline{\P} , -t ) = L ( \overline{\P_{1}} , -t ) + L ( \overline{\P_{2}} , -t ) - L ( \overline{\P_{3}} , -t ) \\ 
                    &\mbox{}& \qquad = (-1)^{n} L ( \P_{1}^{\circ} , t ) + (-1)^{n} L ( \P_{2}^{\circ} , t )  - (-1)^{n-1} L ( \P_{3}^{\circ} , t ) \\ 
                    &\mbox{}& \qquad = (-1)^{n} L ( \P^{\circ} , t ) \ . \end{eqnarray*} 
\hfill {} $\Box$

From the Ehrhart-Macdonald reciprocity law we can now conclude a generalized version of Theorem \ref{multithm}: 
\begin{theorem}\label{genmultithm} Suppose the rational po\-lytope $ \P $ is homeomorphic to an $n$-manifold. Then
$ i_{\P} ( \t ) $ and $ j_{\P} ( \t ) $ are quasipolynomials in $ \t \in \Z^{ m } $, satisfying 
   \[ i_{\P} ( - \t ) = (-1)^{ n } j_{\P} ( \t ) \ . \]
\end{theorem} 
{\it Proof.} 
It suffices to prove that $ i_{\P} ( \t ) $ and $ j_{\P} ( \t ) $ are quasipolynomials. 
In fact, {\it once we know this}, the statement follows from Theorem \ref{recipr}: 
  \[ i_{\P} ( - \t ) = L \left( { \P^{ (\t) } }^{\circ}  , -1 \right) =  (-1)^{ n }  L \left( \overline{\P^{ (\t) }} , 1 \right) =  (-1)^{ n } j_{\P} ( \t ) \ .  \] 
To show that our lattice point count operators are quasipolynomials, it clearly suffices to prove that 
$ i_{\P} ( \t ) $ and $ j_{\P} ( \t ) $ are quasipolynomials in one of the components of $ \t $, say $ t_{1} $. 
Because we leave only this one component variable, we may also assume that $\P$ is convex. 
We make a similar unimodular transformation (which leaves the lattice invariant) as in \cite{beck}: 
we may assume that the defining inequalities for $\P^{ (\t) }$ are
\begin{eqnarray*}\begin{array}{rllllll} a_{ 11 } x_{ 1 } & & & & & \leq & t_{1} \\
                   a_{ 21 } x_{ 1 } &+& \dots &+& a_{ 2n } x_{ n } & \leq & t_{2} \\
                   & \vdots \\
                   a_{ m,1 } x_{ 1 } &+& \dots &+& a_{ m,n } x_{ n } & \leq & t_{m} \ . \end{array} \end{eqnarray*}
(Actually, we could obtain a lower triangular form.) 
Viewing these inequalities as
\begin{eqnarray*}\begin{array}{rllllll} & & & & x_{ 1 } & \leq & \frac{ t_{ 1 } }{ a_{ 11 }  } \\
                   a_{ 22 } x_{ 2 } &+& \dots &+& a_{ 2n } x_{ n } & \leq & t_{ 2 } - a_{ 21 } x_{ 1 }  \\
                   & \vdots \\
                   a_{ m,2 } x_{ 2 } &+& \dots &+& a_{ m,n } x_{ n } & \leq & t_{ m } - a_{ m,1 } x_{ 1 }  \ , \end{array} \end{eqnarray*}
we can compute the number of lattice points in the interior and closure of $ \P^{ ( \t ) } $ as
   \begin{equation}\label{int} i_{\P} ( \t ) = \sum_{ k = s_{1} }^{ \left[ \frac{ t_{ 1 } - 1 }{ a_{ 11 }  }  \right]  } i_{\Q} \left( t_{ 2 } - a_{ 21 } k , \dots , t_{ m } - a_{ m,1 } k \right)  \end{equation}
and
   \begin{equation}\label{clos} j_{\P} ( \t ) = \sum_{ k = s_{2} }^{ \left[ \frac{ t_{ 1 } }{ a_{ 11 }  }  \right] } j_{\Q} \left( t_{ 2 } - a_{ 21 } k , \dots , t_{ m } - a_{ m,1 } k \right)  \ , \end{equation}
respectively. Here $ s_{1} $ and $ s_{2} $ are rational numbers not depending on $ t_{1} $, and 
the $(n-1)$-dimensional polytope $ \Q^{ (\b) } $ is given by 
  \[ \Q^{ (\b) } = \left\{ \x \in \R^{ n-1 } : \ \B \ \x \leq \b \right\} \ , \] 
where 
  \[ \B = \left( \begin{array}{rcl} a_{ 22 } & \dots & a_{ 2n } \\
                   & \vdots \\
                   a_{ m,2 } & \dots & a_{ m,n } \end{array}  \right) \in M_{ (m-1) \times (n-1) } ( \Z ) \ .  \]
The functions $ i_{\Q} $ and $ j_{\Q} $, over which the summations in (\ref{int}) and (\ref{clos}) range, are 
constant in $ t_{1} $. Thus we only need a weak form of Lemma 4 in \cite{beck} 
to deduce that $ i_{\P} ( \t ) $ and $ j_{\P} ( \t ) $ are quasipolynomials in $ t_{1} $. 
\hfill {} $\Box$

At this point, we find it appropriate to remark why we did not simply start the notion of vector-dilated 
polyotopes with this proof, assuming classical Ehrhart-Macdonald reciprocity. The point of \cite{beck} 
(or at least half of it) was really to give an {\it elementary} proof of Theorem \ref{recipr}.
It is for this reason that we chose to build our proof of Theorem
\ref{genmultithm} upon the work in \cite{beck}. The course of the proof looks like the following 
diagram: 
 \[ \mbox{ (\ref{gi}) } \stackrel{ \mbox{ \cite{beck} } }{ \Longrightarrow } \mbox{ Theorem \ref{multithm} } \Longrightarrow \mbox{ Theorem \ref{recipr} } \Longrightarrow \mbox{ Theorem \ref{genmultithm} } \] 


\section{Extending Stanley's theorem}
We conclude by proving the appropriate generalization of Theorem \ref{stan}, essentially in the same way 
Stanley deduced Theorem \ref{stan} from Theorem \ref{recipr}. 
\begin{corollary} Suppose the rational po\-lytope $ \P $ is
homeomorphic to an $n$-manifold. Denote the set of all (closed) facets of $\P$ by $F$, and let $T$ be a subset 
of $F$, such that $ \bigcup_{ \F \in T } \F $ is homeomorphic to an $(n-1)$-manifold. Let 
  \[ j_{\P, T} ( \t ) = \# \left( \left( \P^{ (\t) } - \bigcup_{ \F \in T } \F^{ (\t) } \right) \cap \Z^{ n } \right) \] 
and 
  \[ i_{\P, T} ( \t ) = \# \left( \left( \P^{ (\t) } - \bigcup_{ \F \in F-T } \F^{ (\t) } \right) \cap \Z^{ n } \right) \ . \] 
Then 
   \[ i_{\P, T} ( - \t ) = (-1)^{ n } j_{\P, T} ( \t ) \ . \]
\end{corollary} 
Again, note that Theorem \ref{genmultithm} is the special case $ T = \emptyset $ of this corollary. 

{\it Proof.} By definition, 
  \[ j_{\P, T} ( \t ) = j_{\P} ( \t ) - \sum_{ \F \in T } j_{\F} ( \t )  \] 
and 
  \[ i_{\P, T} ( \t ) = j_{\P} ( \t ) - \sum_{ \F \in F-T } j_{\F} ( \t ) = i_{\P} ( \t ) + \sum_{ \F \in T } i_{\F} ( \t ) \ . \] 
Hence by Theorem \ref{genmultithm}, 
  \[ i_{\P, T} ( - \t ) = (-1)^{n} j_{\P} ( \t ) + \sum_{ \F \in T } (-1)^{n-1} j_{\F} ( \t ) = (-1)^{n} j_{\P, T} ( \t ) \ . \] 
\hfill {} $\Box$


\vspace{5mm}

{\bf Acknowledgements}. I am grateful to Richard Stanley for an eye-opening conversation, which eventually 
lead to the generalization of Theorem \ref{multithm}. 
This paper was written while I was a research associate at Binghamton
University; I would like to thank the Mathematics department at Binghamton
University for their hospitality. 

\nocite{*}
\addcontentsline{toc}{subsubsection}{References}
\bibliography{thesis}
\bibliographystyle{alpha}

 \sc Department of Mathematical Sciences\\
 State University of New York\\
 Binghamton, NY 13902-6000\\
 {\tt matthias@math.binghamton.edu}

\end{document}